\documentclass[10pt,a4paper]{amsart}
\usepackage{amsmath,eucal}
\usepackage{amssymb}

\newcounter{prp}
\newcounter{thm}
\newcounter{cor}
\newcounter{lem}
\newcounter{exl}

\newtheorem{lem}[lem]{Lemma}
\newtheorem{prp}[prp]{Proposition}
\newtheorem{thm}[thm]{Theorem}
\newtheorem{cor}[cor]{Corollary}

\newcommand{\tJ}{\tilde{J}}
\newcommand{\tg}{\tilde{g}}
\newcommand{\dbla}{\langle\!\langle}
\newcommand{\dbra}{\rangle\!\rangle}
\newcommand{\tnb}{\tilde{\nabla}}
\newcommand{\pnb}{\nabla^\perp}

\newcommand{\tR}{\tilde{R}}
\newcommand{\tK}{\tilde{K}}

\newcommand{\nbo}{\nabla^\perp}
\newcommand{\xv}{\xi^v}
\newcommand{\ev}{\eta^v}
\newcommand{\zv}{\zeta^v}

\newcommand{\pR}{R^\perp}

\newcommand{\Xh}{X^h}
\newcommand{\Yh}{Y^h}
\newcommand{\Zh}{Z^h}
\newcommand{\la}{\langle}
\newcommand{\ra}{\rangle}
\newcommand{\rap}{\rangle^{\perp}}

\newcommand{\ot}{\omega_{\theta}}
\newcommand{\oqt}{\omega_{\sqrt{q} \theta}}

\title{ On the geometry of the normal bundle with a metric of Cheeger-Gromoll type}
\author{Wojciech  Koz\l owski}
\begin{document}
\baselineskip=16pt
\maketitle

\begin{abstract}
We investigate the geometry of a normal bundle equipped with a $(p,q)$-metric, i.e., Riemannian metric of Cheeger-Gromoll type, 
to the submanifold of a Riemannian manifold. We derive all natural object as the Levi-Civita connection, curvature tensor, 
sectional and scalar curvature. We prove that under some natural conditions the sectional curvature of this bundle may be 
bounded from below by given arbitrary large
positive constant. Next we investigate $(p,q)$-metrics from the complex geometry point of view. We show when the normal 
bundle can by equipped with a structure of almost Hermitian, almost K\"ahlerian, conformally almost K\"ahlerian or K\"ahlerian manifold.  

\end{abstract}

\section{Introduction and preliminaries} 
The best known example of a natural metric
on the tangent $TM$ to the Riemannian manifold $(M,g)$ has been constructed by S. Sasaki in the late of 1950s.
The geometry of $TM$ with the Sasaki metric $\tg_S$  is proved to be very rigid. For example, O. Kowalski \cite{Kw} showed that $(TM, \tg)$ is locally symmetric iff $(M,g)$ is flat. Moreover, in \cite{MT} E. Musso and F. Tricerri proved that $(TM, \tg_S)$ is of constant scalar curvature iff $(M,g)$ is locally Euclidean. 

Next,
A. A. Borisenko and  A. L. Yampolski showed the full analogy between 
the geometry of the tangent bundle and the geometry of the normal bundle to the submanifold of a Riemannian manifold.

Another example of a natural metric on $TM$  has been given by J. Cheeger and D. Gromoll in \cite{CG}.
This metric $\tg_{Ch-G}$ is known in literature as
 Cheeger-Gromoll metric.  Later the geometry of $(TM, \tg_{Ch-G})$ has been investigated
by M. Sekizawa \cite{Se} and  by S. Gudmundsson and  E. Kappos \cite{GK}. In \cite{GK} it is proved that if the sectional curvature $K$ of $(M,g)$ is constant then  $(TM, \tg_{Ch-G})$ is of positive scalar curvature iff $K$ is bounded by  some constants. Notice also that M. T. K. Abbassi and 
M. Sarih invesigated in \cite{AS} Cheeger-Gromoll metric in the context of Kiling vector fields.

Recently,  the geometry of $TM$ with a very general deformation of the Cheeger-Gromoll metric $\tg_{a,b}$ has been investigated by M. I. Munteanu (\cite{M}). The author of \cite{M} determined the geometry of \((TM,\tg_{a,b})\) and described some aspects of the geometry of the unit tangent bundle $T_1M\subset TM$. 
He also investigated $(TM,\tg_{a,b})$ from
the complex geometry point of view. Applying the method developed by M. Anastasiei in \cite{A}, Munteanu showed 
$TM$ can be equipped with an almost complex structure $\tilde J$ compatible with $\tg_{a,b}$. He proved when \( (TM,g_{a,b},\tilde J)\) is almost Hermitian,
almost K\"ahlerian,   locally conformal  K\"ahlerian, or K\"ahlerian manifold

Independently, in \cite{BLW} M. Benyounes, E. Loubeau and C. M. Wood introduced a some class of natural Riemannian metrics on vector bundles of Cheeger-Gromoll type.
These metrics, ${ h}_{p,q}$, $p,q\in \mathbb{R}$, $q\ge 0$, called {\em $(p,q)$-metrics}  generalize both the
Sasaki metric and the Cheeger-Gromoll metric, but are less general than the metrics introduced by Munteanu.

 $(p,q)$-metrics have been discovered together with some new harmonics maps, but, as the authors showed, the geometry of
$(TM,{ h}_{p,q})$ is of the independent interest \cite{BLW1}. 

Although, metrics from \cite{M} are much more general than $(p,q)$-metrics, $(p,q)$-metrics have very nice geometrical properties. Among results from 
\cite{BLW1} the following seems to be very important:
{\em Suppose that $M$ is a space of constant curvature. Then there exist $p$ and $q$ such that the scalar curvature of $(TM, h_{p,q})$ is strictly positive.}


In the light of \cite{BLW1} and \cite{BY} there appears a question on the geometry of  the the normal bundle with $(p,q)$-metric.
 In this paper we  investigate $(p,q)$-geometry  of the normal bundle $T^\perp L$ to a submanifold $L$ of Riemannian manifold $M$.
Having $(T^\perp L, h_{p,q})$ we derive all natural geometric objects: the Levi-Civita connection $\tnb$, curvature tensor $\tR$, 
sectional curvature $\tK$ 
and scalar curvature $\tilde S$ of 
$h_{p,q}$. We show that this objects are determined by the Levi-Civita connection $\nabla$ of $L$, 
normal connection $\pnb$, normal curvature
 tensor $\pR$ and its adjoint $\hat R$, and the
parameters $p$ and $q$. 

We prove that (Theorem \ref{Thm_Flat}): {\em for every submanifold  $L$ of codimension $\geq 2$, \( (T^\perp L, h_{p,q}) \) is flat iff $p=q=0$, $L$ is flat 
and the normal connection $\pnb$ is flat.} 

We also obtain some estimations for scalar curvature and prove that 
(Theorem \ref{Thm_E_Scal}): {\em if $L$ is of codimension $\geq 2$, its scalar
 curvature and normal curvature tensor are bounden by some positive constant then the scalar curvature of  $T^\perp L$ 
can by bounded below by arbitrary large constant.}
 
Next we study $(p,q)$-metrics from the complex geometry point of view. We prove that under some natural conditions 
it is always possible to introduce almost complex structure $\tilde J$ on $T^\perp L$ compatible with the given $h_{p,q}$, i.e., $(T^\perp L, h_{p,q}, \tilde J)$ is an almost Hermitian manifold (Proposition \ref{Natural_J}). We also prove when $T^\perp L$ is locally conformal almost K\"ahlerian, almost K\"ahlerian (Theorem 
\ref{Almost_Khl}) or K\"ahlerian manifold (Theorem \ref{Thm_Khl}).

\subsection{Our setting}
Our important data are: a Riemannian manifold $(M,g)$ of dimension $\geq 2$, its submanifold $L\subset M$ ( $1\leq \dim L < \dim M$) and the normal bundle
$\pi:T^\perp L\to L$. $g$ induces a Riemannian metric $\la,\ra$ on $L$ and a fibre metric
$\la,\rap$ in $T^\perp L$. Let $\nabla$ be the Levi-Civita connection of
$\la,\ra$, and $\pnb$ be the normal connection in $T^\perp L$ induced from the Levi-Civita connection of $g$.

All maps, vector fields, sections etc. are assumed to be smooth.

We will denote by $X,Y,Z$ vectors / vector fields tangent to $L$, and by $\xi,\zeta,\eta$ members / sections of $T^\perp L$.
Recall that the connection map related to $\pnb$ it is a bundle morphism   $K: T (T^\perp L) \to T^\perp L$, uniquely determined  by
the conditions 
\begin{itemize}
\item[(K1)] For every $x\in L$ and $\theta\in T^\perp L$, $K:T_\theta (T^\perp _x L\to T^\perp_x L$ is a canonical isomorphism.
\item[(K2)] For every section $\xi $ of $T^\perp L$ and vector field $X$ on $L$,  $K(\xi_\ast X) = \pnb_X \xi$.
\end{itemize}
Under the action of $\pnb$, the bundle $T(T^\perp L)$ splits as: 
$T(T^\perp L) = \mathcal H \oplus \mathcal V$, where $\mathcal V = \ker \pi_\ast$ is the {\em vertical bundle} and 
$\mathcal H = \ker K$ is the {\em horizontal bundle}. Every $A\in T_\theta(T^\perp L)$ can be uniquely written as a sum $A=\mathcal H A + \mathcal V A$,
where $\mathcal H A$ and $\mathcal V A$ is the horizontal and vertical part of $A$, respectively.

We denote by $X^h$ and $\eta^v$ a unique horizontal and vertical lift of $X$ and $\eta$, respectively.
Moreover, $\Theta$ denotes the canonical vertical vector field, i.e., for every $\theta\in T^\perp L$, $\Theta_\theta$ is the vertical lift of the vector $\theta$ to the point $\theta$, or equivalently, $\Theta$ is a section of $\mathcal V$ such that $K(\Theta_\theta) = \theta$.

\begin{lem}[\cite{BY}]\label{LemmaBY} Let $\theta\in T^\perp L$. For every vector fields $X,Y$ on $L$ and every sections $\xi,eta$ of $T^\perp L$ we
have
\begin{eqnarray*}
 [\xv,\ev] &=& 0,\quad [X^h,\ev]_\theta = (\pnb_X \eta)^v_\theta,\\
 {} [X^h,Y^h]_\theta &=&  \big([X, Y]\big)^h_\theta - \big(\pR(X,Y)\theta\big)^v_\theta.
 \end{eqnarray*}
\end{lem}

\begin{lem}\label{UL_lemma2} For every vector field $X$ on $L$ and every sections $\xi$ of $T^\perp L$ we have
\begin{eqnarray}
 [\xv.\Theta] &= &\xv,\label{UL_21}\\
 {[}\Xh, \Theta] &=& 0.\label{UL_22}
\end{eqnarray}
\end{lem}
\begin{proof} 
Take $x\in L$ Suppose $\xi_i$, $i=1,\dots, d$, $d={\rm codim \,} L$ is a local frame of $T^\perp L$ in the neighbourhood $U$ of $x$ such that
$\nbo \xi_i=0$ at $x$. We may write
$\Theta = \tilde \theta^i \xi_i^v $ where $\tilde \theta^i (\eta) = \eta^i$ if $\eta = \eta^i \xi_i$, $\eta\in \pi^{-1}(U)$.

\eqref{UL_21}: Suppose that $\xi =  \xi^i\xi_i$. Observe that $\xv \tilde\theta^i = \xi\circ \pi$.
Then we have
\[ [\xv.\Theta] = (\xv \tilde\theta^i) \xi_i^v + \tilde \theta^i [\xv,\xi_i^v]= (\xv \tilde\theta^i) \xi_i^v = (\xi^i\circ \pi)\xi_i^v = \xv. \]

\eqref{UL_22}: Take $\theta\in T^\perp_x L$. Observe that $(\Xh \tilde \theta^i )(\theta)=0$. 
At the point $\theta$ we have  
\[ {[}\Xh, \Theta] = (\Xh \tilde \theta^i)\xi_i^v +
 \tilde \theta^i 
(\nabla^\perp_X \xi_i)^v = 0.\]
Notice that the assumption that $\nbo \xi_i=0$ at $x$ has been used only in the proof of  formula  \eqref{UL_22}.
\end{proof}

We denote by $R$ and $\pR$ the curvature tensor of $\nabla$ and $\pnb$. Let $\hat R$ be adjoint to $\pR$, i.e.,
\[ \la \hat R (\xi,\eta) X, Y \ra = \la \pR (X,Y) \xi, \eta \rap.\]
Notice that by the Ricci equation it follows that
\[ \hat R (\xi,\eta) = - \hat R (\eta, \xi ).\]
As a direct consequence of the definition of $\hat R$ and the fact that $\hat R$ is skew-symmetric we obtain that for every $X,Y,Z,W\in T_x L$ and every
$\xi,\eta,\zeta,\theta\in T^\perp_x L$,
\begin{eqnarray}
\la \hat R(\pR(X,Y)\theta,\xi)Z,W\ra &=& -\la \pR(Z,W)\xi,\pR(X,Y)\theta\rap,\label{hRpR}\\
\la \hat R(\theta,\eta)\hat R(\xi,\zeta)X,Y\ra &=& - \la\hat R(\theta,\eta)Y,\hat R(\xi,\zeta)X\ra\label{hRhR}
\end{eqnarray}
We will also work with the covariant derivatives  $\pnb_X \pR$ and
$\nabla_X \hat R$ defined in a usual way:
\begin{eqnarray*}
\big(\pnb_X \pR\big)(Y,Z)\eta &=& \pnb_X\big( \pR(Y,Z)\eta\big) - \pR(\nabla_X Y,Z)\eta -\pR(Y,\nabla_X Z)\eta - \pR(Y,Z)\pnb_X\eta,\\
\big(\nabla_X \hat R\big)(\xi,\eta)Y &=&
\nabla_X\big( \hat R(\xi,\eta)Y\big) - \hat R(\pnb_X \xi,\eta)Y -\hat R(\xi,\pnb_X \eta)Y
 - \hat R(\xi,\eta)\nabla_XY.
\end{eqnarray*}
Notice that $\pnb_X \pR$ and
$\nabla_X \hat R$ are related as follows:
\[ \la (\nabla_X \hat R )(\xi,\eta)Z,Y\ra = \la (\pnb_X \pR)(Z,Y)\xi,\eta\rap.\]
Moreover, $\pR$ satisfies the 
{\em second Bianchi identity}:
\begin{equation}\label{sBi} 
\big(\pnb_X \pR\big)(Y,Z) +  \big(\pnb_Z \pR\big)(X,Y) + \big(\pnb_Y \pR\big)(Z,X) = 0.
\end{equation}

We will also use the following convention: 
If $\psi$ is a map and $u$ belongs to the domain of $\psi$ we will often write $\psi_u$ instead of $\psi(u)$.
If  $\xi $ is a section of $T^\perp L$. We will often identify $\xi$ with a
 map $T^\perp L\to T^\perp L$,
$\theta \mapsto 
\xi_{\pi(\theta)}$. 
In particular for every sections $\xi$ and $\eta$ we often identify its fibre product $\la \xi,\eta\rap $ with a map
$\theta \mapsto \la \xi,\eta\rap_{\pi(\theta)}$. If $\theta$ is a member of $T^\perp L$, $\la \xi,\theta\rap$ means
$\la \xi_{\pi(\theta)},\theta\rap$.

\subsection{$(p,q)$-metric and its Levi-Civita connection}
Suppose $\dbla,\dbra$ is a natural metric on $T^\perp L$, i.e., $\dbla,\dbra$ is  determined by $g$, and $\pi:(T^\perp L, \dbla,\dbra)
\to (M, \la,\ra)$ is a Riemannian submersion whose horizontal distribution is the kernel of connection map. Directly by Lemma \ref{LemmaBY} and the definition of the Levi-Civita connection it follows:
\begin{lem} Let $\tnb$ be the Levi-Civita connection of  $\dbla,\dbra$, and $\theta \in T^\perp L$ then
\begin{eqnarray}
  \dbla \tnb_{\Xh} \Yh, \Zh \dbra & = &  \langle \nabla_X Y, Z \rangle, \label{Nhhh}\\
2   \dbla \tnb_{\Xh} \ev, \Zh \dbra_\theta & = & \dbla (\pR(X,Z)\theta)^v,\ev \dbra_\theta, \label{Nhvh}\\
2   \dbla \tnb_{\Xh} \Yh, \zv \dbra_\theta & = &  - \dbla (\pR (X, Y) \theta)^v ,  \zv \dbra_\theta,\label{Nhhv}\\
2   \dbla \tnb_{\Xh} \ev, \zv \dbra & =& \Xh \dbla \ev, \zv \dbra + \dbla (\pnb_X \eta)^v, \zv \dbra - \dbla (\pnb_X \zeta)^v,\ev \dbra,\label{Nhvv}\\
2   \dbla \tnb_{\ev} \Yh, \Zh \dbra_\theta & =& \dbla ((\pR(Y,Z)\theta)^v,\ev \dbra_\theta,\label{Nvhh}\\
2   \dbla \tnb_{\ev} \Yh, \zv \dbra & =&  \Yh \dbla \ev,\zv\dbra - \dbla (\pnb_Y\eta)^v,\zv\dbra - \dbla (\pnb_Y \zeta)^v,\ev\dbra,\label{Nvhv}\\
2   \dbla \tnb_{\ev} \zv, \Zh \dbra & =&  -\Zh \dbla \ev,\zv\dbra + \dbla (\pnb_Z \zeta)^v,\ev \dbra + \dbla (\pnb_Z \eta)^v,\zv\dbra,\label{Nvvh}\\
2   \dbla \tnb_{\xv} \ev, \zv \dbra & =& \xv \dbla \ev,\zv \dbra+ \ev \dbla \xv,\zv \dbra  - \zv \dbla \xv,\ev \dbra.\label{Nvvv}
\end{eqnarray}
\end{lem}

The generalized Cheeger-Gromoll metric on $T^\perp L$ it is natural metric $h_{p,q}$ defined as follows:
\[ h_{p,q}(A,B) = \la\pi_\ast A,\pi_\ast B\ra + \omega^p_\theta (\la KA,KB\rap + q\la KA,\theta\rap\la KB,\theta\rap ),\]
where $p\in \mathbb{R}$ and  $q\ge 0$ are fixed and the function $\omega$ is given by
$$\omega_\theta = \frac{1}{1 + \la\theta,\theta\rap}.$$ 
The metric $h_{p,q}$ is also called $(p,q)$-metric.
Notice that $h_{0,0}$ is the Sasaki metric $h_S$ and $h_{1,1}$ is the Cheeger-Gromoll metric $h_{Ch-G}$.\\

\noindent {\sf Remark} Munteanu in \cite{M} considered natural metrics on tangent bundle whose analogue for normal bundle wolud be as follows: 
\[ \tg_{a,b}(A,B) = a(t)\la \pi_\ast A,\pi_\ast B\ra + b(t) \la KA,\theta\rap\la KB,\theta\rap,\quad A,B\it T_\theta\big(T^\perp L\big),\]  
where $a,b\in C^\infty(\mathbb{R}_+$, $a>0$, $b\geq 0$ and $t=t_\theta = \frac 12 \la \theta,\theta \rap$. Obviously, if
$a(t)=\big(1+2t\big)^{-p}$ and 
$b(t)=q\big(1+2t)^{-p}$ then $\tg_{a,b}=h_{p,q}$.\\

For the future considerations, it is convenient to define the function $\omega_{\sqrt q}$
 by $\omega_{\sqrt q}(\theta) = \omega(\sqrt q \theta)$, i.e.,
$$ \oqt = \frac {1}{ 1 +q \la \theta,\theta\rap}.$$

\begin{lem}\label{Lemma_techniczny} Let $\theta \in T^\perp L$ and $x= \pi \theta$. We have
\begin{enumerate}
\item $\Xh\dbla \ev,\xv\dbra = \dbla (\pnb_X \eta)^v, \xv\dbra + \dbla \ev, (\pnb_X \xi)^v \dbra$,
\item $\dbla (\pR(X,Y)\theta)^v,\xi^v \dbra_\theta = \omega^p_\theta \la \hat R(\theta,\xi) X,Y\ra_x$,
\item $ \xv \dbla \ev,\zv\dbra = -2p \omega^{p+1}_\theta \la \xi,\theta\rap \Big(\la \eta,\zeta\rap + q \la \eta,\theta\rap\la \zeta,\theta\rap)
+q\omega^p_\theta (\la \xi,\eta\rap\la \zeta,\theta\rap\\ + \la\xi,\zeta\rap\la\eta,\theta\rap\Big)$
\end{enumerate}

\end{lem}
\begin{proof}
(1) Take a $\pnb$-parallel vector field $e$ such that $e(0) = \theta$ and $\dot e (0) = X^h(\theta)$, and let $\gamma = \pi e$.
Then the derivative of $t\mapsto \omega^p_{e(t)}$ is equal to zero.  Moreover,
\[ \dbla \ev,\xv\dbra_{e(t)} = \omega^p_{e(t)}\big(\la \eta(\gamma(t)),\xi(\gamma(t))\rap +
 q \la \eta(\gamma(t)),e(t)\rap \la \xi(\gamma(t)),e(t)\rap\big).  \]
Consequently,
\begin{eqnarray*}
X^h(\theta) \dbla \ev,\xv\dbra  &=& \left. \frac{d}{dt}\right|_{t=0}   \dbla \ev,\xv\dbra_{e(t)}\\
&=& \left. \frac{d}{dt}\right|_{t=0} \omega^p_{e(t)}\Big(\la \eta(\gamma(t)),\xi(\gamma(t))\rap +
q \la \eta(\gamma(t)),e(t)\rap \la \xi(\gamma(t)),e(t)\rap\Big)\\
&=& \omega^p_{\theta} \Big( \la \pnb_{\dot\gamma(0)} \eta,\xi(x)\rap +  \la \eta(x),  \pnb_{\dot\gamma(0)} \xi\rap\\
&&+ q( \la \pnb_{\dot\gamma(0)} \eta, e(0) \rap \la \xi(x),e(0)\rap+  \la \eta(x),e(0)\rap \la \pnb_{\dot\gamma(0)} \xi, e(0) \rap)\Big)\\
&=& \omega^p_\theta \big(\la \pnb_{X(0)} \eta, \xi(x)\rap + q \la \pnb_{X(0)} \eta, \theta\rap \la \xi(x),\theta\rap \big)\\
&&+ \omega^p_\theta \big(\la \pnb_{X(0)} \xi, \eta(x)\rap + q \la \pnb_{X(0)} \xi, \theta\rap \la \eta(x),\theta\rap\big)\\
&=& \dbla (\pnb_X \eta)^v, \xv\dbra_\theta + \dbla \ev, (\pnb_X \xi)^v \dbra_\theta.
\end{eqnarray*}
(2) Since  $\la \pR(X,Y)\theta,\theta \rap $,
we have
\begin{eqnarray*}\dbla (\pR(X,Y)\theta)^v,\xi^v \dbra_\theta  &=& \omega^p_\theta \big(\la \pR(X,Y)\theta,\xi(x) \rap 
+ q\la \pR(X,Y)\theta,\theta \rap \la \xi(x),\theta\rap\big)\\
&=& \omega^p_\theta \la \pR(X,Y)\theta,\xi(x) \rap\\
&=& \omega^p_\theta \la \hat R(\theta,\xi)X,Y \ra_x .
\end{eqnarray*}
(3) Take any curve $\gamma$ in the fibre $T^\perp_x L$ such that $\gamma(0) = \theta$ and $\dot \gamma (0) = \xv(\theta)$.
Then $\la K\ev(\gamma(t)),K\zv(\gamma(t)\rap = \la\eta(x),\zeta(x)\rap$ and $\la K\ev(\gamma(t)),\gamma(t)\rap = \la \eta(x),\gamma(t)\rap$.
Thus we have
\begin{eqnarray*}
\xv(\theta) \dbla  \ev,\zv\dbra & =& \left. \frac{d}{dt}\right|_{t=0} \dbla  \ev,\zv\dbra_{\gamma(t)}\\
&=&\left. \frac{d}{dt}\right|_{t=0} \Big( \omega^p_{\gamma(t)} (\la\eta(x),\zeta(x)\rap + q\la \eta(x),\gamma(t)\rap\la \zeta(x),\gamma(t)\rap)\Big)\\
&=& -2p \omega^{p+1}_\theta \la \xi,\theta\rap \big(\la \eta,\zeta\rap + q \la \eta,\theta\rap\la \zeta,\theta\rap\big)\\
&&
+q\omega^p_\theta \big(\la \xi,\eta\rap\la \zeta,\theta\rap + \la\xi,\zeta\rap\la\eta,\theta\rap\big).
\end{eqnarray*}
\end{proof}

\begin{lem}\label{LCpq} Let $\theta\in T^\perp L$. The Levi-Civita connection $\tnb$ of the Cheeger-Gromoll type metric $h_{p,q}=\dbla,\dbra$ 
is given by
\begin{eqnarray}
\big( \tnb_{Xh} \Yh\big)_\theta  &=& (\nabla_X Y)^h_\theta - \frac{1}{2}(\pR(X,Y)\theta)^v_\theta,\label{LChh}\\
\big( \tnb_{\Xh} \ev \big)_\theta &=&  \frac{1}{2}\omega^p_\theta (\hat R (\theta,\eta) X)^h_\theta + (\pnb_X\eta)^v_\theta,\label{LChv}\\
\big( \tnb_{\ev} \Yh\big)_\theta & =& \frac{1}{2}\omega^p_\theta (\hat R (\theta,\eta) Y )^h_\theta,\label{LCvh}\\
\big(\tnb_{\xv} \ev\big)_\theta & =& -p\omega_\theta\big(\la \xi,\theta \rap \ev + \la \eta,\theta\rap \xv\big)_\theta\label{LCvv}\\
&&+ \big( \nu_\theta \la \xi,\theta \rap \la \eta,\theta\rap +  \mu_\theta \la \xi,\eta\rap\big) \Theta_\theta,\nonumber
\end{eqnarray}
where $\nu$ and $\mu$ are functions defined by
\begin{eqnarray*} 
\nu_\theta &=& \,\frac{pq\omega_\theta}{q\la \theta,\theta \rap +1},\\
 \mu_\theta &=&\frac{q+p\omega_\theta}{q\la \theta,\theta \rap +1}.
\end{eqnarray*}
\end{lem}
\begin{proof} \eqref{LChh}: It follows directly from \eqref{Nhhh} and \eqref{Nhhv}.

\eqref{LChv}: The equality of the vertical  parts follows from Lemma \ref{Lemma_techniczny} (1) and \eqref{Nhvv}, the equality of the horizontal  parts follows from Lemma \ref{Lemma_techniczny} (2) and \eqref{Nhvh}.

\eqref{LCvh}: The fact that the vertical part of the left hand side is degenerated follows from Lemma \ref{Lemma_techniczny} (1) and \eqref{Nvhv}.
 Thus the equality  follows from Lemma \ref{Lemma_techniczny} (2) and \eqref{Nvhh}.

 \eqref{LCvv}: The horizontal part of the left hand side is degenerated by \eqref{Nvvh} and Lemma \ref{Lemma_techniczny} (1).

 By \eqref{Nvvv} and Lemma \ref{Lemma_techniczny} (3) we get
 \begin{eqnarray}
 \dbla \tnb_{\xv} \ev, \zv \dbra_\theta & =&  \frac{1}{2} \Big(\xv \dbla \ev,\zv \dbra+ \ev \dbla \xv,\zv \dbra  -
\zv \dbla \xv,\ev \dbra\Big)_\theta\label{wzor_fi0}\\
 &=& \omega^p_\theta(q + p\omega_\theta) \la \xi,\eta\rap\la \zeta,\theta\rap \nonumber \\
&&- p\omega^{p+1}_\theta
 (\la \xi,\theta\rap\la \eta,\zeta\rap + \la\eta,\theta\rap\la\xi,\zeta\rap)\nonumber \\
 &&-pq\omega^{p+1}_\theta \la \xi,\theta
\rap\la \eta,
\theta\rap\la\zeta,\theta\rap\nonumber.
 \end{eqnarray}
If $\theta = 0$ then right hand sides of \eqref{wzor_fi0} and \eqref{LCvv} are equal to zero, so our assertion follows.

Suppose that $\theta\ne 0$.
We write $\varphi = (\tnb_{\xv} \ev)_\theta$, for simplicity.
Applying the definition of $\dbla,\dbra$ to the left hand side of  \eqref{wzor_fi0} we
obtain
\begin{eqnarray}
\omega^p_\theta (K\varphi + q \la K\varphi, \theta\rap \theta ) &=& - p\omega^{p+1}_\theta (\la \xi,\theta\rap \eta
+ \la \eta,\theta\rap \xi)\label{wzor_fi}\\
&&-pq\omega^p_\theta \la \xi,\theta\rap\la\eta,\theta\rap \theta\nonumber\\
&& + \omega_\theta^p (p\omega_\theta + q)\la \xi,\eta \rap \theta\label{wzor_fi2}\nonumber.
\end{eqnarray}
Moreover, observe that
\begin{eqnarray*}
K\varphi  &=& \frac{1}{\la \theta,\theta\rap} \la K\varphi , \theta \rap \theta +\Big(K\varphi - \frac{1}{\la \theta,\theta\rap} 
\la K\varphi , \theta \rap \theta\Big),\\
\la \xi,\theta\rap \eta
+ \la \eta,\theta\rap \xi &=& \frac{2}{\la \theta,\theta\rap} \la \xi,\theta\rap \la \eta,\theta\rap\theta\\
&&+ \Big(\la \xi,\theta\rap \eta
+ \la \eta,\theta\rap \xi - \frac{2}{\la \theta,\theta\rap} \la \xi,\theta\rap \la \eta,\theta\rap\theta\Big),
\end{eqnarray*}
where the first terms of the right hand sides of the identities is orthogonal to $\theta$ with respect to the fibre metric $\la,\rap$.
Substituting these identities to \eqref{wzor_fi} we see that the coefficient at $\theta$ must be equal to zero. Consequently,
\[  \la K\varphi ,\theta\rap = \frac{\la \theta,\theta\rap (p\omega_\theta +q)}{q \la \theta,\theta\rap +1}\la \xi,\eta\rap
-\frac{p\omega_\theta(q \la \theta,\theta\rap+2)}{q \la \theta,\theta\rap +1}\la \xi,\theta\rap\la\eta,\theta\rap.\]
Substituting this identity to \eqref{wzor_fi} we obtain \eqref{LCvv}.
\end{proof}

Notice that by \eqref{LCvv} it follows that the fibres $T^\perp_x L$, $x\in L$, are totally geodesic.

\begin{cor}\label{LCcor} For every vector fields $X,Y,Z$ on $L$, every sections $\xi,\eta,\zeta$ of $T^\perp L$, and every $\theta\in T^\perp L$ we have
\begin{eqnarray*}
\big(\tnb_{[\xv,\zv]} \Zh\big)_\theta &=& 0,\\
\big(\tnb_{[\xv,\zv]} \ev\big)_\theta &=& 0,\\
\big(\tnb_{[\Xh,\Yh]} \Zh\big)_\theta &=& \big(\nabla_{[X,Y]} Z
- \frac{\omega_\theta^p}{2}\hat R(\theta,\pR(X,Y)\theta)Z\big)^h_\theta\\
&&-\frac 12 \big(\pR([X,Y],Z)\theta\big)^v_\theta,\nonumber\\
\big(\tnb_{[\Xh,\Yh]} \ev\big)_\theta &=& \frac 12 \omega^p_\theta \big( \hat R(\theta,\eta)[X,Y]\big)^h_\theta - \mu_\theta \la \pR(X,Y)\theta,\eta\rap \Theta_\theta\\
&&+\big( \pnb_{[X,Y]}\eta + p\omega_\theta \la \eta,\theta\rap \pR(X,Y)\theta \big)^v_\theta,\nonumber\\
\big(\tnb_{[\Xh,\xv]} \Zh \big)_\theta&=& \frac 12 \omega^p_\theta \big(\hat R(\theta,\pnb_X \xi) Z\big)^h_\theta,\\
\big(\tnb_{[\Xh,\xv]} \ev\big)_\theta &=& -p\omega_\theta\big(\la \pnb_X\xi,\theta\rap \eta + \la\eta,\theta\rap\pnb_X\xi\big)^v_\theta\\
&&+ \big(\nu_\theta \la \pnb_X\xi,\theta\rap\la \eta,\theta\rap + \mu_\theta\la \pnb_X \xi,\eta\rap\big)_\theta \Theta_\theta
\end{eqnarray*}
\end{cor}

\begin{proof} Direct consequence of Lemmas \ref{LCpq} and \ref{LemmaBY}.  
\end{proof}

\subsection{Kowalski type  lemma}
If $\pi': E' \to N$ and $\pi'':E''\to N$ are  vector bundles then
a (smooth) fibre preserving ($\pi'' B = \pi'$)  map $B: E'\to E''$ is called a {\em bundle morphism} if for every $x\in N$, $B: E'_x\to E''_x$ is 
linear.
The following is a version of the Kowalski's Lemma \cite[p. 125]{Kw}
for $T^\perp L$. Since the proof of Lemma \ref{KTLemma} can be based on the same method as Kowalski's one, the demonstration is omitted.

Suppose that
$F:T^\perp \to T^\perp L$ and $G: T^\perp L \to TL$ are bundle morphisms then we define the vertical lift $F^v$ and the horizontal lift $G^h$
of $F$ and $G$ as follows:
\[ F^v(\theta) = \big[F(\theta)\big]^v_\theta,\quad G^h(\theta) =\big[ G(\theta)\big]^h_\theta. \]
Notice that the canonical vertical vector field $\Theta$ is nothing but the vertical lift of the identity map ${\rm Id}: T^\perp L \to T^\perp L$.

\begin{lem}\label{KTLemma} Let $\theta\in T^\perp L$ and $\pi\theta =x$. Then for every vector $\xi\in T^\perp_x L$, we have
\begin{eqnarray*}
\big(\tnb_{\xv} F^v\big)_\theta & =&
\big(F(\xi)\big)^v_\theta -
 p\omega_\theta\big(\la\xi,
\theta \rap \big(F(\theta)\big)^v_\theta
+\la F(\theta),\theta\rap \xv_\theta \big)\\
&&+\mu_\theta\la \xi,
F(\theta)\rap \Theta_\theta + \nu_\theta
\la \xi,\theta\rap  \la F(\theta),\theta\rap \Theta_\theta\nonumber\\
\big(\tnb_{\xv} G^h\big)_\theta & =&
\big(G(\xi)\big)^h_\theta +\frac{\omega^p_\theta}{2}\big(\hat R(
\theta,\xi)
G(\theta)\big)^h_\theta.
\end{eqnarray*}
Let $X\in T_x L$. If $\eta$ satisfies $(\pnb_X\eta) = 0$ and $\eta_x = \theta$ then
\begin{equation*}
\big(\tnb_{\Xh}
F^v\big)_\theta =
\big(\tnb_{\Xh} (F\circ\eta)^v \big)_\theta,\quad \big(\tnb_{\Xh} G^h\big)_\theta = \big(\tnb_{\Xh} (G\circ\eta)^h \big)_\theta.
\end{equation*}
\end{lem}

\begin{proof}
We will prove above formulae for $F$. Proof  for $G$  is analogous. Let $A\in T_\theta (T^\perp L)$. 
Denote by $d$ the  codimension of $L$.
In a neighbourhood $U$ of $x$ we may write
\[ F = \sum_{i=1}^d \lambda^i \xi_i,\]
where $(\xi_i)$ is a local orthonormal frame of $T^\perp U$, 
and $\lambda_i\in C^\infty(\pi^{-1}(U))$.
Since $F$ is a bundle morphism, $\lambda_i$ restricted to the fibre $T^\perp_x L$ is a linear functional . Notice that $\lambda_i = \la F, \xi_i\rap$.
We have 
\[ \tnb_A F^v  =
 \sum_{i=1}^d (A\lambda^i) 
(\xi_i)^v_\theta +\sum_{i=1}^d \lambda^i(\theta) \tnb_A \xi_i^v.\]
Let $\alpha$ be a curve in $T^\perp_x L$ defined by $\alpha(t)= \theta + t \xi$. Then $\xi^v_\theta = \dot \alpha(0)$.
Thus taking $A= \xi^v_\theta$ we get
$A\lambda_i=\lambda^i(\xi)$.
Consequently, 
\[\big(\tnb_{\xv} F^v\big)_\theta = F^v_\xi + \sum_{i=1}^d \la F,\xi_i\rap_\theta  \big( \tnb_\xi^v \xi_i^v\big)_\theta.\]
Now the formula follows by Lemma \ref{LCpq}.

Next, let  $\eta$ be a section of $T^\perp L$ such that $\eta(x)=\theta$ and $(\pnb_X \eta)_x = 0$. Then $X^v_\theta = \eta_\ast X$.
Thus, taking $A=X^v_\theta$ we get
$A \lambda^i = (\eta_\ast X)\lambda^i = X(\lambda^i\circ \eta)$ and $\lambda^i(\theta = (\lambda^i\circ \eta)(x)$.
Consequently, 
\[ \big(\tnb_{\Xh} F^v \big)_\theta
 = \sum_{i=1}^d X(\lambda^i
\circ 
\eta)(\xi_i)^v_\theta +
 \sum_{i=1}^d(\lambda^i\circ \eta)(x) \big(\tnb_{X^h} 
\xi_i^v\big)_\theta\]
On the other hand, $(F\circ \eta)^v_\theta =  (\lambda^i\circ \eta\circ \pi)(\theta)(\xi_i)^v_\theta$ and  
$X^h_\theta(\lambda^i\circ \eta\circ \pi)  = X(\lambda^i\circ \eta)$. Therefore,
\[ \big(\tnb_{\Xh} (F\circ \eta)^v\big)_\theta = \sum_{i=1}^d X^h_\theta(\lambda^i\circ \eta\circ \pi) (\xi_i)^v_\theta + 
\sum_{i=1}^d(\lambda^i\circ\eta\circ\pi)(\theta) \big(\tnb_{\Xh} \xi_i^v\big)_\theta.\]  
Hence the second formula is proved.
\end{proof}

\noindent {\sf Remark} In the case of the tangent bundle of a Riemannian manifold equipped with the Sasaki metric above lemma is due to 
O. Kowalski \cite{Kw}. However, his proof is based on some other method than ours.\\

Let $Y,Z$ be vector fields on $L$ and $\zeta$ be a section of $T^\perp L$.
Consider bundle morphisms ${\rm Id}:T^\perp L\to T^\perp L$, $\pR(Y,Z):T^\perp L\to T^\perp L$ and  $\hat R(\zeta)(Y):T^\perp L\to TL$, where
$\hat R(\zeta)(Y)\theta = \hat R(\zeta,\theta)Y$. Then by Lemma \ref{KTLemma} we have
\begin{eqnarray}
\big(\tnb_{\xv} \Theta\big)_\theta &=& (1-p\omega_\theta \la \theta,\theta \rap )\xv_\theta+q \oqt \la \xi,\theta\rap \Theta_\theta,\label{KTLvT}\\
\big(\tnb_{\xv} \pR(Y,Z)^v\big)_\theta &= & \big(\pR(Y,Z)\xi)^v_\theta
-p\omega_\theta\la \xi,\theta\rap\big(\pR(Y,Z)\theta\big)^v_\theta\label{KTLvpR}\\
&&+\mu_\theta \la \xi,\pR(Y,Z)\theta\rap \Theta_\theta\nonumber\\
\big(\tnb_{\xv} \hat R(\zeta)(Y)^h\big)_\theta &=& \big(\hat R(\zeta,\xi)Y\big)^h_\theta +\frac{\omega^p_\theta}{2}
\big(\hat R(\theta,\xi)\hat R(\zeta,\theta)Y\big)^h_\theta,\label{KTLvhR}\\
\big(\tnb_{\Xh} \Theta\big)_\theta &=& 0,\label{KTLhT}\\
\big( \tnb_{\Xh} \pR(Y,Z)^v\big)_\theta&=& \frac{\omega_\theta^p}{2}\big(
\hat R(\theta,\pR(Y,Z)\theta)X\big)^h_\theta\label{KTLhpR}\\
&&+\big((\pnb_X\pR)(Y,Z)\theta\big)^v_\theta\nonumber\\
&&+\big(\pR(\nabla_X Y,Z)\theta+\pR(Y,\nabla_X Z)\theta\big)^v_\theta,\nonumber\\
\big(\tnb_{\Xh} \hat R(\zeta)(Y)^h\big)_\theta &=& \frac 12\big(\pR(X,\hat R(\theta,\zeta)Y)\theta\big)^v_\theta\label{KTLhhR}\\
&&+\big((\nabla_X\hat R)(\zeta,\theta)Y+\hat R(\pnb_X \zeta,\theta)Y+\hat R(\zeta,\theta)\nabla_X Y\big)^h_\theta.\nonumber
\end{eqnarray}

\section{The geometry of $T^\perp L$}
\setcounter{lem}{0}
\setcounter{thm}{0}
\setcounter{prp}{0}
\setcounter{equation}{0}
\setcounter{cor}{0}
\setcounter{clm}{0}

In this section we assume that $T^\perp L$ is equipped with the  Cheeger-Gromoll type metric $h_{p,q}=\dbla,\dbra$.

\subsection{Curvature tensor and sectional curvature}

\begin{lem}\label{tRKK} For every vector fields $X,Y,Z$ on $L$, every sections $\xi,\eta,\zeta$ of $T^\perp L$ and every $\theta\in T^\perp L$ we have
\begin{eqnarray*}
\big(\tnb_{\Xh}\tnb_{Yh} \Zh\big)_\theta &=& \big(\nabla_X \nabla_Y Z\big)^h_\theta 
-\frac{\omega_\theta^p}{4}\big(
\hat R(\theta,\pR(Y,Z)\theta)X\big)^h_\theta\\
&&-\frac 12 \big((\pnb_X\pR)(Y,Z)\theta+\pR(\nabla_X Y,Z)\theta\\
&&+\pR(Y,\nabla_X Z)\theta
+\pR(X,\nabla_Y Z)\theta\big)^v_\theta,\\
\big(\tnb_{\Xh}\tnb_{\Yh} \zv\big)_\theta &=& \big(\pnb_X \pnb_Y \zeta \big)^v_\theta 
-\frac{\omega^p_\theta}{4}\big(\pR(X,\hat R(\theta,\zeta)Y)\theta\big)^v_\theta\\
&&-\frac{\omega^p_\theta}{2} \big((\nabla_X\hat R)(\zeta,\theta)Y+\hat R(\pnb_X \zeta,\theta)Y\\
&&+\hat R(\zeta,\theta)\nabla_X Y
+\hat R(\pnb_Y \zeta,\theta) X\big)^h_\theta\\
\big(\tnb_{\Xh}\tnb_{\zv} \Yh\big)_\theta &=& 
\frac{\omega^p_\theta}{4}\big(\pR(\hat R(\theta,\zeta)Y,X)\theta\big)^v_\theta\\
&&-\frac{\omega^p_\theta}{2}\big((\nabla_X\hat R)(\zeta,\theta)Y+\hat R(\pnb_X \zeta,\theta)Y+\hat R(\zeta,\theta)\nabla_X Y\big)^h_\theta\\
\big(\tnb_{\xv}\tnb_{\Yh} \Zh \big)_\theta &=& \frac{\omega^p_\theta}{2}\big( \hat R(\theta,\xi)\nabla_Y Z\big)^h_\theta\\
&&-\frac 12 \big(\pR(Y,Z)\xi)^v_\theta
-p\omega_\theta\la \xi,\theta\rap\big(\pR(Y,Z)\theta\big)^v_\theta\\
&&-\frac 12 \mu_\theta \la \xi,\pR(Y,Z)\theta\rap \Theta_\theta\\
\big(\tnb_{\xv}\tnb_{\ev} \Zh \big)_\theta &=& 
 \frac{\omega^p_\theta}{2} \big( \hat R (\xi,\eta) Z
+ \frac{\ot^p}{2} \hat R(\theta,\xi)\hat R(\theta,\eta)Z-2p\omega_\theta \la \xi, \theta \rap \hat R(\theta,\eta) Z\big)^h_\theta,\\
 \big(\tnb_{\xv}\tnb_{\Zh} \ev \big)_\theta &=& \frac{\omega^p_\theta}{2} \big( \hat R (\xi,\eta) Z
+   \frac{\ot^p}{2} \hat R(\theta,\xi)\hat R(\theta,\eta)Z-2p\omega_\theta \la \xi, \theta \rap \hat R(\theta,\eta) Z\big)^h_\theta\\
&&- p\omega_\theta \big(\la \xi,\theta\rap \pnb_Z \eta + \la \pnb_Z \eta,\theta\rap \xi\big)^v_\theta\\
&&+\big(\nu_\theta \la\xi,\theta\rap \la \pnb_Z \eta,\theta\rap +\mu_\theta \la \xi,\pnb_Z \eta\rap \big) \Theta_\theta,\\
 \big(\tnb_{\Xh}\tnb_{\xv} \ev \big)_\theta &=& \frac{p\omega_\theta^{p+1}}{2}\big( \la \xi, \theta\rap \hat R(\eta,\theta) X +\la \eta,\theta\rap \hat R (\xi,\theta) X\big)^h_\theta\\
 &&-p\omega_\theta \big( \la \pnb_X \xi,\theta\rap \eta +\la \pnb_X \eta,\theta\rap \xi+\la \xi,\theta\rap \pnb_X \eta +\la \eta,\theta\rap \pnb_X \xi\big)^v_\theta\\
 &&+\nu_\theta \big(
\pnb_X \xi,
\theta\rap \la \eta,\theta\rap +\la \xi,\theta\rap \la \pnb_X \eta,\theta\rap\big) \Theta_\theta\\
&&+\mu_\theta\big(\la \pnb_X\xi,\eta\rap + \la \pnb_X \eta,\xi\rap\big)\Theta_\theta,\\
\big(\tnb_{\zv}\tnb_{\xv} \ev \big)_\theta &=& -p \ot \big( \la \zeta,\xi\rap - \ot \la \zeta,\theta\rap \la \xi,\theta\rap (p+2)\big)\ev\\
&&-p \ot \big( \la \zeta,\eta\rap - \ot \la \zeta,\theta\rap \la \eta,\theta\rap (p+2)\big)\xv\\
&&+\big(p\ot(p\ot+\mu_\theta)\la \xi,\theta\rap\la\eta,\theta\rap-
(p(1- \ot)-1)\mu_\theta
\la\xi,\eta\rap\big)\zv\\
&&-\oqt\big(\nu_\theta(3-2\ot)+\oqt(q^2+2p\ot^2))\la \xi,\eta\rap \la \zeta,\theta\rap \Theta_\theta\\
&&-p^2\ot^2\oqt\big(\la \eta,\theta\rap \la \zeta,\xi\rap + \la\xi,\theta\rap\la\zeta,\eta\rap\big)\Theta_\theta\\
&&-\nu_\theta\big(2\nu_\theta +3\mu_\theta+\ot\oqt(2-p-2q)\big)\la\xi,\theta\rap\la\eta,\theta\rap\la\zeta,\theta\rap\Theta_\theta
\end{eqnarray*}
\end{lem}

\begin{proof}
All formul\ae {} are consequence of Lemma \ref{LCpq} and formulae \eqref{KTLvT}-\eqref{KTLhhR}. However, a proof of the last identity is laborious.

Notice also the symmetry $\tnb_{\zv}\tnb_{\xv}\ev= \tnb_{\zv}\tnb_{\ev}\xv$ that follows from the facts that $[\ev,\xv]=0$ and $\tnb$ is torsion free.
\end{proof}

\begin{prp}\label{Curvature_Tensor_R} Let $X,Y,Z$ be  vector fields on $L$ and $\xi,\eta,\zeta$ be sections of $T^\perp L$.
The curvature tensor $\tR$ of $\tnb$ at the point $\theta\in T^\perp L$ is given by: 
\begin{eqnarray*}
\tR_\theta (X^h,Y^h)Z^h & =& \big(R(X,Y)Z)^h_\theta\\
&&-\frac{\omega^p_\theta}{4}
\Big( \hat R(\theta,\pR(Y,Z)\theta)X-
\hat R(\theta,\pR(X,Z)\theta)Y -
2\hat R(\theta,
\pR(X,Y)\theta)Z\Big)^h_\theta\nonumber \\
 &&+ \frac{1}{2}\Big(\big(\pnb_ZR\big)(X,Y)\theta\Big)^v_\theta,\nonumber\\
 \tR_\theta (X^h,Y^h)\ev &=& \frac{\omega^p_\theta}{2} \Big(\big(\nabla_X \hat R\big)(\theta,\eta)Y - \big(\nabla_Y \hat R
\big)(\theta,\eta)X\Big)^h_\theta\\
 &&+\big(\pR (X,Y)\eta\big)^v_\theta  +\frac{\omega^p_\theta}{4}\Big( \pR(Y,\hat R (\theta,\eta)X)\theta - \pR(X,\hat R(\theta,\eta)Y)\theta\Big)^v_\theta\\
 &&-p\omega^p_\theta \la \eta,\theta\rap \big(\pR(X,Y)\theta\big)^v_\theta +\mu_\theta\la \pR(X,Y)\theta,\eta\rap
\Theta_\theta,\\
 \tR_\theta (X^h,\ev)Z^h &=& \frac{\omega^p_\theta}{2}\big( (\nabla_X \hat R)(\theta,\eta)Z\big)^h_\theta +
\frac{1}{2}\big( \pR(X,Z)\eta \big)^v_\theta\\
 &&-\frac{p\omega_\theta}{2} \la \eta,\theta\rap \big( \pR(X,Z)\theta\big)^v_\theta - \frac{\omega^p_\theta}{4}\big(\pR(X,\hat R(\theta,\eta)Z)\theta\big)^v_\theta\\
 &&+\frac{1}{2}\mu_\theta\la \pR(X,Z)\theta,\eta\rap \Theta_\theta\\
\tR_\theta (X^h, \ev)\xv & =& \frac{p\omega^{p+1}_\theta}{2}\Big(\la \eta,\theta\rap
\hat R (\theta,\xi) X - \la\xi,\theta\rap  \hat R(\theta,\eta)X\Big)^h_\theta  -\frac{\omega^p_\theta}{2} \big( \hat R(\eta,\xi)X\big)^h_\theta\\
&&-\frac{\omega^{2p}_\theta}{4}\big( \hat R (\theta,\eta)\hat R(\theta,\xi)X\big)^h_\theta,\\
\tR_\theta (\xv,\ev)Z^h & =& \omega^p_\theta\big(\hat R(\xi,\eta)Z\big)^h_\theta + p\omega^{p+1}_\theta \Big(\la \eta,\xi\rap \hat R(\theta,\xi)Z -\la 
\xi,\theta\rap
\hat R(\theta,\eta)Z\Big)^h_\theta\\
&&+\frac{\omega^{2p}_\theta}{4}\Big(\hat R (\theta,\xi) \hat R (\theta,\eta)Z - \hat R (\theta, \eta)\hat R(\theta,\xi)Z\Big)^v_\theta,\\
\tR_\theta (\zv,\xv)\ev &=& \big( \nu_\theta-\ot (2\nu_\theta + p(p-2)\oqt \ot \big)\la\eta,\theta\rap \la\zeta,\theta\rap \xv_\theta\\
&&-\oqt \big(q+p\ot(2\ot-(p-2)(1-\ot))\big)\la\zeta,\eta\rap_\theta\xv_\theta\\
&&-\big( \nu_\theta-\ot (2\nu_\theta + p(p-2)\oqt \ot \big)\la\eta,\theta\rap \la\xi,\theta\rap \zv_\theta\\
&&+\oqt \big(q+p\ot(2\ot-(p-2)(1-\ot))\big)\la\xi,\eta\rap_\theta\zv_\theta\\
&&+\oqt\big(\oqt(q^2-p(p-2)\ot^2)+\nu_\theta((p-2)\ot+3-p)\big)\\
&&\cdot\big(\la \zeta,\eta\rap\la\xi,\theta\rap-\la\xi,\eta\rap\la\zeta,\theta\rap\big)_\theta
\Theta_\theta
\end{eqnarray*}
\end{prp}

\begin{proof} All formulae are a direct consequence of Lemma \ref{tRKK} and Corollary \ref{LCcor}. However, in a proof of the first formula we 
also apply the second  Bianchi identity \eqref{sBi}.
\end{proof}
Notice, that $\tR (\xv, \ev)\zv$ can be written in the form
\begin{eqnarray}
\tR_\theta (\xv, \ev)\zv &=& a_\theta \la \zeta,\theta\rap \big(\la \eta, \theta\rap \xv - \la \xi,\theta\rap \ev\big)_\theta\label{tRabc}\\
&&+ b_\theta \big(\la \eta, \zeta\rap \xv - \la \xi,\zeta\rap \ev\big)_\theta\nonumber\\
&&+c_\theta \big( \la \xi,\theta\rap \la \eta,\zeta\rap - \la \eta , \theta \rap \la \xi , \zeta \rap \big)_\theta\Theta_\theta,\nonumber
\end{eqnarray}
where $a_\theta$, $b_\theta$ and $c_\theta$ are given by
\begin{eqnarray*}
a_\theta & =&  \frac{p\omega^2_\theta(p+q-2 - q\la \theta,\theta\rap)}{1+q\la \theta,\theta\rap},\\
b_\theta  &= & \frac{2p\omega_\theta - p^2\la \theta,\theta\rap \omega^2_\theta + q}{1 + q\la\theta,\theta\rap },\\
c_\theta &= &\frac{pq\omega_\theta -q^2 +\omega_\theta^2\big(p^2 - 2p(1+q) + pq (p-4)\la\theta,\theta\rap\big)}{\big(1 + q\la\theta,\theta\rap\big)^2}.
\end{eqnarray*}
Observe that these functions satisfies the identity
\begin{equation}\label{Eqabc}
a_\theta-q b_\theta = \oqt c_\theta.
\end{equation}

\begin{prp}\label{Curvature_K} Given $X,Y\in T_x L$ and $\xi,\eta\in T^\perp_x L$ be pairs of orthogonal vectors of unit length.
Let $\theta \in T^\perp_x L$. 
 Then the sectional curvature $\tK$ of $\dbla,\dbra$ at $\theta$ is given by:
\begin{eqnarray}
\tK(\Xh_\theta\wedge \Yh_\theta) & =& K(X\wedge Y) - \frac{3}{4}\omega^p_\theta |\pR(X,Y)\theta|_\perp^2,\label{tKXY}\\
\tK(\Xh_\theta\wedge \ev_\theta) & =& \frac{1}{4}\frac{\omega^p_\theta}{1 + (q \la\eta,\theta\rap)^2}|\hat R(\theta,\eta)X|^2,\label{tKXe}\\
\tK (\xv_\theta\wedge \ev_\theta) &=&
\frac{1}{\ot^{p}}\frac{b_\theta + a_\theta\big[(\la \xi,\theta\rap)^2 + 
(\la \eta,\theta\rap)^2\big] }{1+q\big[(\la \xi,\theta\rap)^2 + (\la \eta,\theta\rap)^2\big]}\label{tKxe},
\end{eqnarray}
where $K$ denotes the sectional curvature of $\la,\ra$,  $|\,\,|_\perp^2=\la,\rap$ 
and $|\,\,|^2=\la,\ra$. Note that if $L$ is dimension one (resp. codimension one) submanifold then
$\tK(\Xh\wedge \Yh)$ (resp. $\tK (\xv\wedge \ev)$) is omitted. Moreover, in these cases \(\tK(\Xh_\theta\wedge \ev_\theta) = 0  \).
\end{prp}

\begin{proof} Write, for simplicity,  $\xv$, $\ev$, $\Xh$ and $\Yh$ instead of $\xv_\theta$, $\ev_\theta$, $\Xh_\theta$ and $\Yh_\theta$. 

The formulae \eqref{tKXY}- \eqref{tKxe} are a consequence of 
Proposition \ref{Curvature_Tensor_R} and the identities \eqref{hRpR} and \eqref{hRhR}. However, 
applying O'Neill's formulae (see \cite{Be}, Chapter 9 \S D) we may prove \eqref{tKXY} and \eqref{tKXe} in much easier way.

Since the fibres of $T^\perp L$ are totally geodesic the O'Nill's formulae becomes:
\begin{eqnarray*}
\tK(\Xh\wedge \Yh)&=&K(X\wedge Y) - \frac{3}{4} \| \mathcal V \big([\Xh,\Yh]\big)\|^2,\\
\tK(\Xh\wedge \ev) &=& \frac{1}{\|\ev\|^2} \Big( \| \mathcal H \big(\tnb_{\Xh} \ev\big)\|^2 - \| \mathcal V\big( \tnb_{\ev} \Xh\big)\|^2\Big).
\end{eqnarray*} 
Since, $\| \mathcal V \big([\Xh,\Yh]\big)\|^2 = \ot^p|\pR(X,Y)\theta|^2_\perp$ formula \eqref{tKXY} follows.
Since,
$\tnb_{\ev} \Xh$ is horizontal, $\| \mathcal H \big(\tnb_{\Xh} \ev\big)\|^2 = \frac 14 \omega^{2p} |\hat R(\theta,\eta)X|^2$ and
$\|\ev\|^2 =\ot^p( 1+\big(\la \eta,\theta\rap\big)^2)$, the formula \eqref{tKXe} follows.

We prove formula \eqref{tKxe}. We have
\begin{equation*}
\tK(\xv\wedge\ev) = \frac {\dbla \tR (\xv,\ev)\ev,\xv\dbra}{\|\xv\|^2\|\ev\|^2-\dbla \xv,\ev\dbra^2}.
\end{equation*}
Denote by $N$ and $D$ the numerator and the denominator of  \(\tK(\xv\wedge\ev)\), respectively. Applying the definition of $\dbla,\dbra$
we obtain
\begin{eqnarray*}
\|\xv\|^2 &=& \ot^p \big(1+ q \big(\la \xi,\theta\rap\big)^2\big),\\
\|\ev\|^2 &=& \ot^p \big(1+ q \big(\la \eta,\theta\rap\big)^2\big),\\
\dbla \xv,\ev\dbra &=& q\ot \la \xi,\theta\rap\la\eta,\theta\rap,\\
\dbla \xv,\Theta_\theta\dbra &=& \ot^p \la \xi,\theta\rap \big( 1 + q \la \theta,\theta\rap\big).
\end{eqnarray*}
Consequently,
\begin{equation}\label{tKD}
D = \ot^{2p} \big(1+q \big(\big(\la \xi,\theta\rap\big)^2 +\la \eta,\theta\rap\big)^2\big)
\end{equation}

Next, by \eqref{tRabc}, we obtain
\begin{equation*}
\tR(\xv,\ev)\ev = \big( a_\theta \big(\la \eta,\theta\rap\big) + b_\theta\big)\xv - a_\theta\la\xi,\theta\rap\la \eta,\theta\rap \ev + c_\theta\la \xi,\theta\rap \Theta_\theta.
\end{equation*}
It follows that
\begin{equation}\label{tKN}
N=  \ot^p\big( a_\theta \big(
\la\eta,\theta\rap\big)^2+
(qb_\theta +(1+q \la\theta,\theta\rap)c_\theta)\big(\la \xi,\theta\rap\big)^2 + b_\theta\big)
\end{equation}
Now \eqref{tKxe} is a direct consequence of \eqref{tKD}, \eqref{tKN} and \eqref{Eqabc}.
\end{proof}

\begin{cor} Suppose $\dim L=1$ and $\dim M = 2$. Then for arbitrary $p\in \mathbb{R}$ and $q\geq 0$, $(T^\perp L, h_{p,q})$ is flat.
\end{cor}

As a consequence of Proposition \ref{Curvature_Tensor_R} and Proposition \ref{Curvature_K} we obtain (compare \cite[Theorem 1]{BY}):

\begin{thm}\label{Thm_Flat} Suppose that ${\rm codim\,} L \geq 2$. The manifold $(T^\perp L,h_{p,q})$ is flat iff $p=q=0$, $(L,\la,\ra)$ is flat and the normal connection $\pnb$ is flat.
\end{thm}

\begin{proof} ($\Leftarrow$) A straightforward consequence of  Proposition \ref{Curvature_Tensor_R}.

($\Rightarrow$)   We have $0=\tR_\theta (\Xh,\Theta)\Zh $ for every   $X,Y$ and $\theta$. Consequently, for every $X,Z$ and $\theta$
\[ \frac{1}{2}\left( 1- \frac{p}{1+\la\theta,\theta\rap}\right) (\pR ( X,Z)\theta)^v_\theta = 0.\]
It follows that $\pR \equiv 0$, i.e.,  the normal connection $\pnb$ is flat.
Next, for every vector fields $X,Y$ and $Z$ we have $0 = \tR_\theta(\Xh,\Yh)\Zh = (R(X,Y)Z)^h_\theta$. Thus $R\equiv 0$, so $(L,\la,\ra)$ is flat.

Now we must show that $p=q=0$. By Proposition \ref{Curvature_K} it follows that
\begin{equation}\label{Flat_Curvature_Eq}
a_\theta(\la \xi,\theta\rap)^2 + (\la \eta,\theta\rap)^2) +b_\theta = 0,
\end{equation}
for every orthonormal $\xi,\eta$ (${\rm codim\,}L \geq 2$) and every $\theta$.
Applying the definition of $a_\theta$ and $b_\theta$, and substituting $\theta = 0$ to \eqref{Flat_Curvature_Eq} we obtain that $q=-2p$. Keeping this in mind, and substituting $\theta = \xi$  to \eqref{Flat_Curvature_Eq} we get $p=0$, and then $q=0$.
\end{proof}

\subsection{Scalar curvature of $T^\perp L$} Assume that $d=\dim L$ and $d'= {\rm codim\,} L$.
 Take $\theta\in T^\perp_x L$.
First we construct an orthonormal basis of $T_\theta( T^\perp L)$. Consider two cases (i) $\theta = 0$ and (ii) $\theta\ne 0$, separately.

(i) 
Since $\theta=0$ we see that 
$$\dbla A,B\dbra = \la \pi_\ast A,\pi_\ast B\ra + \la KA,KB\rap, \quad A,B\in T_\theta(T^\perp L).$$
Take an orthonormal basis 
$(X_1,\dots,X_d)$ of 
$T_x L$ and an orthonormal basis 
$(\xi_1,\dots,\xi_{d'})$ of $T^\perp_x L$.
Put $E_i= (X_i)^h_\theta$, $i=1,\dots,d$ and $F_i=(\xi_i)^v_\theta$, $i=1,\dots,d'$. then one can see that $(E_1,\dots,E_d,F_1,\dots,F_{d'})$
is an orthonormal basis of $T_\theta(T^\perp L)$.

(ii) If $\theta\ne 0$ we proceed as follows: Take an orthonormal basis 
$(X_1,\dots,X_d)$ of 
$T_x L$ and an orthonormal basis 
$(\xi_1,\dots,\xi_{d'})$ of $T^\perp_x L$ where $\xi_1= (1/|\theta|_\perp) \theta$.
Put $E_i= (X_i)^h_\theta$, $i=1,\dots,d$ and 
\begin{equation*}
\left \{
\begin{array}{lclcc}
F_1 &=& ({\ot^p \oqt})^{-\frac 12} (\xi_1)_\theta^v & & \\
F_i& =& \ot^{-\frac p2} (\xi_i)_\theta^v &\text{for} & i=2,\dots,d'
\end{array}
\right.
\end{equation*}
We show that $(E_1,\dots,E_d,F_1,\dots,F_{d'})$
is an orthonormal basis of $T_\theta(T^\perp L)$.

Obviously, $\dbla E_i,E_j\dbra =
 \delta_{i j}$ and 
$\dbla E_i,
F_j\dbra = 0$.
Write $F_i= \delta_i \ot^{-\frac p2} (\xi_i)_\theta^v $, where $\delta_1= \oqt^{-\frac 12}$ and $\delta_i= 1$ for $i=2,\dots,d'$
\begin{eqnarray*}
\dbla F_i,F_j\dbra &=& \delta_i\delta_j \ot^{-p} \dbla (\xi_i)_\theta^v,(\xi_j)_\theta^v\dbra\\
&=& \delta_i\delta_j \ot^{-p} \ot^p ( \la \xi_i,\xi_j\rap + q \la \xi_i,\theta\rap\la \xi_j,\theta\rap)\\
&=& \delta_i\delta_j(\delta_{ij}+ q \la \theta,\theta\rap \delta_{i1}\delta_{1j}).
\end{eqnarray*}
Now applying the definition of $\delta_i$ we easily check that $\dbla F_i,F_j\dbra = \delta_{ij}$.

\begin{thm}\label{ScalarTL} Adopt above notation.
Denote by $S$ and $\tilde S$ the scalar curvature of $(L,\la,\ra)$ and $(T^\perp L,\dbla,\dbra)$. Take $\theta\in T^\perp_x L$
then
\begin{eqnarray}
\tilde S &=& S_x -  \frac{3}{4}\omega^p_\theta
 \sum_{i,j=1}^d
|\pR(X_i,X_j)\theta|_\perp^2 
+\frac 12 \ot^p \sum_{i=1}^d\sum_{j=1}^{d'} |\hat R(\theta,\xi_j) X_i|^2  \label{tS}\\
&& + (d'-1)\ot^{-p} \oqt \big[ 2a_\theta \la\theta,\theta\rap + b_\theta (d'+(d'-2)q\la\theta,\theta\rap )\big],\nonumber
\end{eqnarray}
where $d$ and $d'$ denote the dimension and codimension of $L$.

Notice that if $\dim L=1$ (resp. $ {\rm codim\,} L = 1$ then the first  (resp. last) term in $\tilde S$ is omitted. Moreover, if $\dim L=1$ and $\dim M=2$ then $\tilde S = 0$.
\end{thm}

\begin{proof}
Clearly it suffices to prove \eqref{tS} for $\theta\ne 0$.

Let $\mathcal E=(E_1,\dots,E_d,F_1,\dots,F_{d'})$ denote an orthonormal basis of $T_\theta(T^\perp L)$ from (ii).
We have
\begin{eqnarray*}
\tilde S &=& 2\sum_{i< j} \tK(E_i\wedge E_j) + 2\sum_{i j} \tK(E_i\wedge F_j) +2\sum_{i<j}\tK(F_i\wedge F_j)=2\sigma_1+2\sigma_2+2\sigma_3.
\end{eqnarray*}
By the definition of $\mathcal E$ it follows that
\begin{eqnarray*}
\sigma_1 &=&  \sum_{i<j} \tK((X_i)^h_\theta,(X_j)^h_\theta),\\
\sigma_2 &=& \sum_{i, j} \tK((X_i)^h_\theta, (\xi_j)^v_\theta),\\
\sigma_3 &=& \sum_{i<j} \tK ((\xi_i)^v_\theta,(\xi_j)^v_\theta).
\end{eqnarray*}
Applying \eqref{tKXY} we obtain
\begin{eqnarray*}
\sigma_1 &=& \sum_{i<j} \Big(
K(X_i\wedge X_j) - \frac{3}{4}\omega^p_\theta |\pR(X_i,X_j)\theta|_\perp^2\Big) \\
&=& \frac 12 \left( S_x -\frac{3}{4}\omega^p_\theta
 \sum_{i,j}
|\pR(X_i,X_j)\theta|_\perp^2 \right).
\end{eqnarray*}
Applying \eqref{tKXe} we obtain
\begin{eqnarray*}
\sigma_2 &=& \frac 14 \ot^p \sum_i\left( \sum_j\frac {1}{1+(q\la \xi_j\theta\rap)^2}
|\hat R(\theta,\xi_j) X_i|^2\right)\\
&=&\frac 14 \ot^p \sum_{i,j}|\hat R(\theta,\xi_j) X_i|^2.  
\end{eqnarray*}
Applying \eqref{tKxe} we obtain
\begin{eqnarray*}
\sigma_3 &=& \sum_{i<j} \frac{\ot^{-p}}{1+q\la \theta,\theta\rap (\delta_{1i}+\delta_{1j})} \big(a_\theta \la\theta,\theta\rap
(\delta_{1i}+\delta_{1j})+b_\theta \big)\\
&=& \sum_{i<j} \frac{\ot^{-p}}{1+q\la \theta,\theta\rap \delta_{1i}} \big(a_\theta \la\theta,\theta\rap
\delta_{1i}+b_\theta \big)\\
&=&  (d'-1)\ot^{-p}\oqt (a_\theta \la\theta,\theta\rap +b_\theta) +\frac{(d'-1)(d'-2)}{2}\ot^{-p} b_\theta\\
&=& \frac 12 (d'-1)\ot^{-p} \oqt \big( 2a_\theta \la\theta,\theta\rap + b_\theta (d'+(d'-2)q\la\theta,\theta\rap )\big).
\end{eqnarray*}
This finishes the proof.
\end{proof}

\subsection{Estimates of the scalar curvature} 

State the step key
\begin{lem}\label{lem_EST_SCAL}
Let $c_1,c_2\in \mathbb{R}$ and $d'\in \mathbb{N}$, $d\geq 2$. There exist $p\in \mathbb{R}$ and $q\geq 0$ such that the following function $\Phi=\Phi_{c_1,c_2}$
is strictly positive for  $t\geq 0$.
\[\Phi(t) = c_1-c_2 \frac{t}{(1+t)^p} +\frac{(1+t)^{p-2}}{(1+qt)^2} P(t).\]
Here $P$ is a polynomial of form
\[ P(t) = \alpha_0 t^3 + \alpha_1 t^2+\alpha_2 t + \alpha_3,\]
where 
\begin{eqnarray*}
\alpha_0 &=& (d'-2) q^2,\\
\alpha_1 &=& -q(6p+2p^2-4q+d'(1+2p+2q-p^2)),\\
\alpha_2 &=& 2q(p+d')+q(d'-2)(2p+q)+p(d'-2)(2-p),\\
\alpha_3 &=& d'(2p+q).
\end{eqnarray*}
\end{lem}

Let $C>0$. We say that $\pR$ (resp. $\hat R$) is bounded by $C$ if for every $x\in  L$, every   $X,Y\in T_x L$ and every $\xi,\eta\in T^\perp_x L$
$|\pR(X,Y)\xi|_\perp \leq C |X||Y||\xi|_\perp$ (resp. $|\hat R(\xi,\eta) X \leq C |\xi|_\perp |\eta |_\perp |X|$).

\begin{thm}\label{Thm_E_Scal} Let $D>0$.
Assume that $L$ is a submanifold of $(M,g)$ codimension $\geq 2$. Suppose that there exists $C>0$ such that $|S|<C$, $|\pR|<C$,
where $S$ denotes the scalar curvature of $L$.
Then there exists a $(p,q)$-metric $h_{p,q}$ on $T^\perp L$ such that the scalar curvature of $(T^\perp L,h_{p,q})$ is $> D$.
\end{thm}

\begin{proof} By theorem \ref{ScalarTL} we have
\begin{eqnarray*}
 \frac 1{d'-1}\big(  \tilde S(\theta)-D\big)  
& >& -\frac 1{d'-1}(C+D) - \frac 3{4(d'-1)} d'^2 C \ot^p |\theta|_\perp^2\\
&&+\ot^{-p} \oqt \big( 2a_\theta \la\theta,\theta\rap + b_\theta (d'+(d'-2)q\la\theta,\theta\rap )\big)\\
&=&\Phi_{c_1,c_2}
( | \theta |_\perp^2),
\end{eqnarray*} 
where $c_1= -\frac 1{d'-1}(C+D) $ and $c_2 = \frac 3{4(d'-1)} d'^2 C$. Now the assertion follows by Lemma \ref{lem_EST_SCAL}. 
\end{proof}

\begin{cor}
Suppose $L$ is compact submanifold of codimension 
$\geq 2$, then there exist 
$p\in \mathbb{R}$ and $q\geq 0$ such that the scalar curvature of $(T^\perp L,h_{p,q})$ is strictly positive.
\end{cor}

\begin{cor}
Suppose $L$ is one-dimensional submanifold, i.e., $L$ is a curve, in $M$ with $\dim M \geq 3$ then then there exist 
$p\in \mathbb{R}$ and $q\geq 0$ such that the scalar curvature of 
$(T^\perp L,h_{p,q})$ 
is strictly positive.
\end{cor}

\section{A natural almost complex structure on $T^\perp L$}
\setcounter{lem}{0}
\setcounter{thm}{0}
\setcounter{prp}{0}
\setcounter{equation}{0}
\setcounter{cor}{0}
\setcounter{clm}{0}

We want to find a  natural almost complex structure on $T^\perp L$ compatible with the given $(p,q)$-metric 
$h_{p,q}=\dbla,\dbra$. Clearly, $\dim T^\perp L$ can be odd, so in general it is impossible. Therefore, we restrict our consideration to the 
following case:

{\em $M$ is $2k$-dimensional manifold equipped with an almost complex structure $J$, compatible with the metric $g$, and $L$ is a $k$-dimensional totally
real submanifold in the sense that: for every $x\in L$, $T_x M$ splits as a direct orthogonal sum: $T_xM = T_xL \oplus^\perp J(T_x L)$.}

\subsection{An almost complex structure}
Modify a method form \cite{M} we will seek an almost complex structure of the form:
\begin{eqnarray*}
(\tJ \xv)_\theta &=& a (J\xi)^h_\theta + b\la \xi,\theta\rap J^h_\theta,\\
(\tJ X^h)_\theta &=& c (JX)^v_\theta + d \la JX,\theta\rap \Theta_\theta,
\end{eqnarray*}
where $a,b,c,d$ are functions on $T^\perp L$, and $J^h$ denotes the horizontal lift of the bundle morphism $J:T^\perp L \to TL$.
Moreover, we may suppose that the functions $a$ and $c$ are non-negative.

We must have $\tJ^2 = -1$ and $h_{p,q}\circ \tJ = h_{p,q}$. Writing each from these equalities for horizontal and vertical vectors we get
\begin{eqnarray*}
ac &=&1,\\
ad + b(c + d \la \theta,\theta\rap) &=& 0,\\
cb + d (a + b\la \theta,\theta\rap) & =& 0,
\end{eqnarray*}
from the first one, and from the second one
\begin{eqnarray*}
a^2 & =& \omega^p,\\
c^2\omega^p &=& 1,\\
a^{-2}(2ab+b^2\la\theta,\theta\rap) &=& q,\\
qd^2 (\la \theta,\theta\rap)^2 + (d^2 + 2cdq)\la \theta,\theta\rap + 2cd + c^2q &=&0.
\end{eqnarray*}
It follows that:
\begin{eqnarray*}
a_\theta  &=& \omega^\frac{p}{2},\quad
c_\theta  =  \omega^{-\frac{p}{2}},
\end{eqnarray*}
\begin{eqnarray*}
b_\theta  &=& -\omega^\frac{p}{2}\frac{1 \pm \sqrt{1 + q \la\theta,\theta\rap}}{\la\theta,\theta\rap},\quad
d_\theta = - \omega^{-\frac{p}{2}}\frac{q\la\theta,\theta\rap + 1 \pm \sqrt{1 + q\la\theta,\theta\rap}}{\la\theta,\theta\rap(1+ q\la\theta,\theta\rap)}.
\end{eqnarray*}
If we choose the minus sign in $b$ or $d$ then we obtain a singularity at $\theta = 0$. Since we want $\tJ$ to be defined on the whole $T^\perp L$ we
choose minus sign in `$\pm$' in both $b$ and $d$. Then we may write on the whole $T^\perp L$:

\begin{eqnarray*}
b_\theta &=&  \omega^\frac{p}{2}\frac{q}{1+\sqrt{1+q\la\theta,\theta\rap}},\quad
d_\theta = -\frac{1}{\omega^\frac{p}{2}} \frac{q}{1+q\la\theta,\theta\rap +\sqrt{1 + q\la\theta,\theta\rap}}.
\end{eqnarray*}

We proved:

\begin{prp}\label{Natural_J}
 $(T^\perp L, \dbla,\dbra, \tJ)$ is an almost Hermitian manifold if $\tJ$ is given by:
\begin{eqnarray*}
(\tJ \xv)_\theta &=&  \omega^\frac{p}{2}_\theta \left( (J\xi)^h_\theta +
\frac{q}{1+\sqrt{1+q\la\theta,\theta\rap}} \la \xi,\theta\rap J^h_\theta\right),\\
(\tJ X^h)_\theta &=&  \frac{1}{\omega^\frac{p}{2}_\theta} \left(
(JX)^v_\theta -
 \frac{q}{1+q\la\theta,\theta\rap +\sqrt{1 + q\la\theta,\theta\rap}}
 \la JX,\theta\rap \Theta_\theta\right).
\end{eqnarray*}
\end{prp}

Notice that if $p=q=0$ then  $\tJ \xv =(J\xi)^h$ and $\tJ \Xh = (JX)^v$.

\subsection{On the integrability of $\tJ$}\label{SS_Integrability_of_J}

In this section we assume additionally  that $(M,J,g)$ is K\"ahlerian and $L$ is totally geodesic.
Moreover, let $\nabla^g$ and $R^g$ denote the Levi-Civita connection and curvature tensor of $g$, and let $K$ be the sectional curvature of
$(L,\la,\ra)$.  Recall that on a K\"ahler manifold, $J$ is
parallel, i.e., $\nabla^g J=0$, and $R^g$ satisfies the following (\cite[Proposition 4.5]{KN}):
\begin{equation}\label{Properties_of_Rg}
R^g(V,W)(JU) = J(R^g(V,W)U),\quad R^g(JV,JW)U=R(V,W)U,
\end{equation}
for every vector fields $V,W,U$ on $M$.

Suppose $\tJ$ is given as in Proposition \ref{Natural_J}. We ask whether $\tJ$ is integrable. Let $\tilde N$ be the torsion of $\tJ$.
Recall, that $\tilde N$ is a tensor field of form:
\[ \frac{1}{2} \tilde N(A,B):= [\tJ A, \tJ B] - [A,B]  - \tJ [A, \tJ B] - \tJ [\tJ A, B].\]
\begin{thm} Suppose $(T^\perp L, \tJ, h_{p,q})$ is Hermitian. Then $L$ is the space of constant curvature
\[ K = 2^{p-1}\frac{p+p\sqrt{1+q} + 2q}{\sqrt{1+q}(1+\sqrt{1+q})}\]
\end{thm}

\begin{proof}
Since $\tilde N (\Xh,\Yh) = 0$, applying Lemma \ref{UL_lemma2} and the identity $\nabla^g J = 0 $, after derivations we get
\[ \pR(X,Y)\theta = \Phi(\theta) \big( \la JY,\theta\rap JX - \la JX,\theta\rap JY\big),\]
where
\[ \Phi(\theta) = \frac{(1+\la\theta,\theta\rap)^{p-1}(p+p\sqrt{1+q\la\theta,\theta\rap} +q + q\la\theta,\theta\rap)}
{\sqrt{1+q\la\theta,\theta\rap}(1 + \sqrt{1 + q\la\theta,\theta\rap})}.\]

Suppose now that $X,Y$ are orthogonal and have unit length. Putting $\theta = JY$ and applying the above we have
\begin{eqnarray*}
\la \pR (X,Y)JY,JX\rap &=& \Phi(JY) \big(\la JY,JY\rap\la JX,JX\rap - (\la JX,JY\rap)^2\big)\\
&=& \Phi(JY).
\end{eqnarray*}
On the other hand applying Ricci equation, \eqref{Properties_of_Rg} and the fact that $L$ is totally geodesic we get
\begin{eqnarray*}
\la \pR (X,Y)JY,JX\rap &=& g( R^g (X,Y)JY,JX)\\
&=& g(J(R^g (X,Y)Y),JX)\\
&=& g( R^g (X,Y)Y,X)\\
&=& K(X \wedge Y).
\end{eqnarray*}
Since $\la JY,JY\rap = \la Y,Y\ra = 1$, the assertion follows.
\end{proof}

\subsection{Fundamental form}
 In the present section we adopt the assumptions from Section \S \ref{SS_Integrability_of_J}. 
 
 Let $\varphi$ be the fundamental form of $(T^\perp L, \tJ, h_{p,q})$, where   $h_{p,q}=\dbla,\dbra$:
 \[ \varphi(A,B) = \dbla A, \tJ B\dbra,\]
 for every $A,B\in T(T^\perp L)$. By the definition of $\tJ$ it follows directly that
 \begin{eqnarray} 
\varphi(\Xh,\Yh) &=& 0,\qquad
\varphi(\xv,\ev) = 0,\label{fi_E1}\\
\varphi_\theta(\Xh,\xv) &=& 
\omega^\frac{p}{2}_\theta \left( \la X, J\xi\ra - 
\frac{q}{1 +\sqrt{1 + q\la\theta,\theta\rap}}
\la \xi,
\theta\rap\la X,J\theta
\ra\right).\label{fi_E2}
\end{eqnarray}

\begin{lem}\label{Lemma_dfi}
We have
\begin{eqnarray}
(d\varphi)(\xv,\ev,\zv) & =& 0,\label{dfivvv}\\
(d\varphi)(\Xh,\Yh,\Zh) & = & 0,\label{dfihhh}\\
(d\varphi)(\Xh,\Yh,\xv) &=& 0,\label{dfihhv}\\
(d\varphi)_\theta(\xv,\ev,\Xh) & =&  \omega^\frac{p}{2}_\theta \left( p\ot - \frac{q}{1+\sqrt{1 + q \la \theta,\theta\rap}}\right)
\big( \la \xi,\theta\rap \la X,J\eta\ra\label{dfivvh}\\ 
&&- \la \eta, \theta\rap \la X, J\xi\ra\big)_\theta\nonumber
\end{eqnarray}
\end{lem}
\begin{proof}
Recall that
for every vector fields $A,B,C$
\begin{eqnarray}
 (d\varphi)(A,B,C) &=& A \varphi(B,C) - B \varphi (A,C) + C \varphi (A,B)\label{Koszul_Formula}\\
 &&+\varphi([A,C],B) - \varphi ([A,B],C) - \varphi([B,C],A)\nonumber
 \end{eqnarray}

\eqref{dfivvv}: It follows directly from the second identity in \eqref{fi_E1} and Lemma
\ref{LemmaBY}.

\eqref{dfihhh}: Let $A=\Xh$, $B= \Yh$ and $C= \Zh$. By \eqref{fi_E1} if follows that the three first terms on the right hand side of \eqref{Koszul_Formula} are equal to zero.
Applying Lemma \ref{LemmaBY} and the equality $\la \pR(X,Y)\theta,\theta\rap = 0$
we derive that
\begin{eqnarray*}
\varphi_\theta([A,B],C) &=& \varphi_\theta([X,Y]^v,\Zh)\\
&=& \varphi_\theta(\Zh, (\pR(X,Y)\theta)^v)\\
&=&\omega^\frac{p}{2}_\theta \la Z, J\pR(X,Y)\theta\ra 
\end{eqnarray*}
Since $(M,g,J)$ is K\"ahlerian and $L$ is totally geodesic, by Ricci Equation it follows that 
\[ \la Z, J\pR(X,Y)\theta\ra = - \la R^g(X,Y)Z,J\theta\ra.\]
Consequently,
\[  (d\varphi)_\theta (A,B,C) = \omega^\frac{p}{2}_\theta \la R^g(X,Y)Z +R^g(Z,X)Y + R^g(Y,Z)X,J\theta\ra = 0,\] 
by the first Bianchi identity for $R^g$.

\eqref{dfihhv}: Keeping in mind that $\nabla^g J = 0$ one can obtain
\begin{eqnarray*}
\Xh_\theta \varphi(\Yh,\xv) &=& 
\omega^\frac{p}{2}_\theta 
(\la 
\nabla_X Y,J\xi\ra 
+ \la Y, J\pnb_X \xi\ra)_\theta \\
&&- \omega^\frac{p}{2}_\theta  \frac{q}{1 + \sqrt{1 + q \la\theta\,\theta\rap}} (\la \pnb_X \xi,\theta\rap\la Y,J\theta\ra + \la \xi,\theta\rap \la 
\nabla_X Y,J\theta\ra),
\end{eqnarray*}
\begin{eqnarray*}
\Yh_\theta \varphi(\Xh,\xv) &=& 
\omega^\frac{p}{2}_\theta 
(\la 
\nabla_Y X,J\xi\ra 
+ \la X, J\pnb_Y \xi\ra)_\theta \\
&&- \omega^\frac{p}{2}_\theta  \frac{q}{1 + \sqrt{1 + q \la\theta\,\theta\rap}} (\la \pnb_Y \xi,\theta\rap\la X,J\theta\ra + \la \xi,\theta\rap \la 
\nabla_Y X,J\theta\ra)_\theta,
\end{eqnarray*}
\begin{eqnarray*}
\varphi_\theta ([\Xh,\Yh],\xv) &=& \omega^\frac{p}{2}_\theta \left( \la[X,Y],J\xi\ra_\theta -   \frac{q}{1 + \sqrt{1 + q \la\theta\,\theta\rap}}
\la \xi,\theta\rap\la[X,Y],J\theta\ra\right),\\
\varphi_\theta ([\Xh,\xv],\Yh) &=& -\omega^\frac{p}{2}_\theta \left( \la Y,J\pnb_X\xi\ra_\theta -  \frac{q}{1 + \sqrt{1 + q \la\theta\,\theta\rap}} \la\pnb_X \xi,\theta\rap \la Y, J\theta\ra\right),\\
\varphi_\theta ([\Yh,\xv],\Xh) &=& -\omega^\frac{p}{2}_\theta \left( \la X,J\pnb_Y\xi\ra_\theta -  \frac{q}{1 + \sqrt{1 + q \la\theta\,\theta\rap}} \la\pnb_Y \xi,\theta\rap \la X, J\theta\ra\right)
\end{eqnarray*}
Now by these and \eqref{Koszul_Formula}, \eqref{dfihhv} follows.

\eqref{dfivvh}: Directly by Lemma \ref{LemmaBY}, \eqref{fi_E1} and  \eqref{Koszul_Formula} it follows that
\[ (d\varphi)(\xv,\ev,\Xh) = \xv \varphi(\ev, \Xh) - \ev \varphi(\xv, \Xh).\]
Now \eqref{dfivvh} is a straightforward consequence of simply differentiation.
\end{proof}
Let $\alpha$ be a one-form on $T^\perp L$ determined by the conditions:
\begin{eqnarray*}
\alpha(\Xh) &= &0,\\
\alpha_\theta (\xv) &=& -\left( p\omega_\theta -  \frac{q}{1 + \sqrt{1 + q \la\theta\,\theta\rap}}\right)\la \xi,\theta\rap.
\end{eqnarray*}

\begin{thm}\label{Almost_Khl}
The almost Hermitian manifold $(T^\perp L, \tJ, h_{p,q})$ is locally conformal almost K\"ahlerian, i.e.,
\[d\varphi = \alpha \wedge \varphi.\]
In particular, $(T^\perp L, \tJ, h_{p,q})$ is almost K\"ahlerian iff $p=q=0$.
\end{thm}
\begin{proof}
The assertion follows directly from Lemma \ref{Lemma_dfi}, the definition of $\alpha$ and the identity
\[ (\alpha\wedge \varphi)(A,B,C) = \alpha(A)\varphi(B,C) - \alpha(B)\varphi(A,C) + \alpha(C)\varphi(A,B),\]
for every $A,B,C\in \ T(T^\perp L)$.
\end{proof}

\begin{lem}\label{N_for_p=0_q=0}
Suppose $p=0$ and $q=0$ then 
\begin{eqnarray*}
\tilde N_\theta (\Xh,\Yh) &=& 2(\pR (X,Y)\theta)^v_\theta,\\
\tilde N_\theta (\xv,\ev )&=& 2(\pR (J\eta,J\xi)\theta)^v_\theta,\\
\tilde N_\theta (\Xh,\xv )&=& 2(J\pR (X,J\xi)\theta)^h_\theta.
\end{eqnarray*}
\end{lem}

\begin{thm}\label{Thm_Khl}
$(T^\perp L, \tJ, h_{p,q})$  is K\"ahlerian iff $p=0$, $q=0$ and $\pnb $ is flat.
\end{thm}

\begin{proof} 
A direct consequence of Lemma \ref{N_for_p=0_q=0} and Theorem \ref{Almost_Khl}.
\end{proof}

\vskip20pt
\noindent{\sc Wojciech Koz\l owski}\\
Faculty of Mathematics and Computer Science, University of \L\' od\' z\\
ul. Banacha 22, 90-238 \L\' od\' z, Poland\\
e-mail: {\tt wojciech@math.uni.lodz.pl}
\vskip10pt

\end{document}